\documentclass[conference]{ieeetran}
\usepackage{cite}
\usepackage[pdftex]{graphicx}
\usepackage{epstopdf}
\usepackage{amsmath}
\usepackage{cases}
\usepackage[caption=false,font=footnotesize]{subfig}
\usepackage{fixltx2e}
\usepackage{amssymb}
\usepackage{mathtools}
\usepackage{bm}
\usepackage{multirow}
\usepackage{algorithm}
\usepackage{algorithmic}
\usepackage{amsthm}
\usepackage{mathrsfs}
\usepackage{color}
\usepackage{textcomp,mathcomp}

\begin{document}

\title{Coordinated Active-Reactive Power Management of ReP2H Systems with Multiple Electrolyzers}

\author{
  \IEEEauthorblockN{Yangjun~Zeng\textsuperscript{1}, Buxiang~Zhou\textsuperscript{1}, Jie~Zhu\textsuperscript{1}, Jiarong~Li\textsuperscript{2}, Bosen~Yang\textsuperscript{2}, Jin~Lin\textsuperscript{2}, and Yiwei~Qiu\textsuperscript{1*}}
  \IEEEauthorblockA{
  \textsuperscript{1}College of Electrical Engineering, Sichuan University, Chengdu 610065, China\\
  \textsuperscript{2}Department of Electrical Engineering, Tsinghua University, Beijing 100084, China\\
  \textsuperscript{*}Email: ywqiu@scu.edu.cn
  }
}

\maketitle

\begin{abstract}
  Utility-scale renewable power-to-hydrogen (ReP2H) production typically uses thyristor rectifiers (TRs) to supply power to multiple electrolyzers (ELZs). They exhibit a nonlinear and non-decouplable relation between active and reactive power. The on-off scheduling and load allocation of multiple ELZs simultaneously impact energy conversion efficiency and AC-side active and reactive power flow. Improper scheduling may result in excessive reactive power demand, causing voltage violations and increased network losses, compromising safety and economy. To address these challenges, this paper first explores trade-offs between the efficiency and the reactive load of the electrolyzers. Subsequently, we propose a coordinated approach for scheduling the active and reactive power in the ReP2H system. A mixed-integer second-order cone programming (MISOCP) is established to jointly optimize active and reactive power by coordinating the ELZs, renewable energy sources, energy storage (ES), and var compensations. Case studies demonstrate that the proposed method reduces losses by 3.06\% in an off-grid ReP2H system while increasing hydrogen production by 5.27\% in average.
\end{abstract}


\section{Introduction}
\label{sec:intro}


Renewable power to hydrogen (ReP2H) is a promising approach for large-scale integration of renewable energy and driving the low-carbon energy transition \cite{yang2022breaking}. Many demonstration projects have been initiated worldwide. Among them, the off-grid ReP2H systems faces fewer grid integration restrictions than the grid-connected setups, offering higher planning and operation flexibility \cite{ibanez2023off}. However, due to the lack of external grid support, it relies on its internal power sources and scheduling of electrolyzers (ELZs) to maintain active and reactive power balance, ensuring safety and economic viability \cite{kafetzis2020energy}.

Due to the limited capacity of individual ELZs, utility-scale ReP2H systems typically consist of multiple ELZs \cite{li2023multi}. The scheduling of multiple ELZs, involving on-off switchings and load allocation\cite{matute2021multi,varela2021modeling,xing2020intermodule,QIU2023EXTEND}, enable flexible adjustment of hydrogen production loads to accommodate varying power generation. For example, \cite{matute2021multi} develops a techno-economic analysis tool for evaluating the cost of P2H considering on-standby-off switchings of the ELZs; \cite{varela2021modeling} proposes a scheduling scheme for multiple alkaline ELZs to address fluctuations in renewable energy and electricity price; \cite{xing2020intermodule} develops intermodule power and thermal management strategies for solid oxide electrolysis cells; while \cite{QIU2023EXTEND} develops a scheduling approach for alkaline ELZs that considers internal mass and heat transfer constraints to exploit the flexibility of ReP2H further.

However, existing research predominantly focuses on active power management, often overlooking the impact of reactive loads from the ELZs. Driven by technological and cost considerations, many ReP2H projects use thyristor rectifiers (TRs) to supply power to the ELZs \cite{ruuskanen2020power}. The TRs draws significant reactive power during hydrogen production, and the non-decoupled relation between active and reactive loads \cite{ruuskanen2020power,koponen2021comparison} may result in very different reactive load as the active load varies. The combination of on-off states and load allocation among the ELZs simultaneously impact the P2H energy conversion efficiency and the power flow in the network.

Specifically, traditional ELZ scheduling strategies \cite{varela2021modeling,QIU2023EXTEND} primarily aim to optimize the energy conversion efficiency. However, we find a conflict between increasing the efficiency and reducing the reactive load of the ELZs, which may lead to excessive reactive power demands, resulting in voltage violations, increased network losses, and unintended negative impacts on safety and the economy. An intuitive example of this conflict is shown in Section \ref{sec:conflict}. While some studies have modeled and analyzed the reactive load of TRs for ELZs \cite{ruuskanen2020power}, they are limited to the rectifier level. Currently, the academic and industrial sectors lack scheduling methods for off-grid ReP2H systems that consider both active and reactive power.

To address these challenges, this paper extends our preliminary works on hydrogen plant scheduling \cite{QIU2023EXTEND,li2023multi} and draws inspiration from the reactive power optimization in distribution networks (DN) \cite{chai2019hierarchical}. We propose a coordinated active-reactive power management method for off-grid ReP2H systems. First, a detailed active and reactive load model of the ELZ is established, taking into account on-off switchings and dynamic thermal effects. Subsequently, we analyze the trade-offs between P2H energy conversion efficiency, voltage, and network losses in scheduling multiple ELZs. Furthermore, leveraging the nonlinear relationship between the active and reactive load of the ELZs, by means of the combination of on-off states and load allocation, and considering available regulating resources such as renewable energy sources (RESs), energy storage (ES), and var compensations, we propose a coordinated active-reactive power management method. This employs a mixed-integer second-order cone programming (MISOCP) to balance energy conversion efficiency and network losses, ultimately enhancing the overall profitability of off-grid ReP2H systems.

The remainder is organized as follows. Section \ref{sec:model} establishes the ReP2H system model and explore the trade-offs between active and reactive power management. Section \ref{sec:management} presents the coordinated scheduling approach. Finally, Section \ref{sec:case} shows the case studies, and Section \ref{sec:conclusion} concludes this work.

\section{Modeling Active and Reactive Loads for Multiple Electrolyzers}
\label{sec:model}

\begin{figure}[tb]
  \centering
  \includegraphics[width=3.3in]{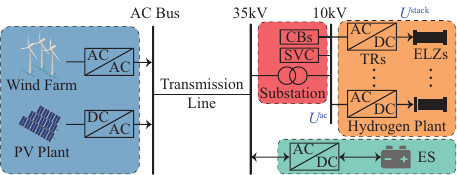}\vspace{-3.5pt}
  \caption{Diagram of a typical off-grid renewable power to hydrogen system.}
  \label{fig:system}\vspace{-10pt}
\end{figure}

The off-grid ReP2H system consists of wind turbines (WTs) and photovoltaics (PVs) plants on the source side, ELZs on the load side, and the ES for balance regulation, as shown in Fig. \ref{fig:system}. Among these, the ES plays a crucial role in maintaining safe operation and participates in production scheduling to support the smooth operation of the hydrogen plant \cite{ibanez2023off}.

The electrolytic load of the ELZ is jointly determined by the DC current and temperature of the stack. This, in turn, determines the AC-side reactive load of the TR. Controlling the active and reactive loads of the ELZ is challenging due to strong coupling and nonlinearity. While WTs, PVs, ES, and compensators supply reactive power, excessive reactive power flow can lead to voltage issues and increased network losses. To explore this, Sections \ref{sec:efficiency} and \ref{sec:real-reactive} investigate the nonlinear efficiency and reactive load characteristics of the ELZ, and Section \ref{sec:conflict} provides a visual depiction of the trade-offs between active and reactive power management.

\subsection{Active Load and Efficiency of Hydrogen Production}
\label{sec:efficiency}

In application, the external characteristics of the electrolysis stack can be described by the empirical U-I curve \cite{ulleberg2003modeling}:
\begin{equation}
\begin{split}
U^{\mathrm{stack}}=& U^{\mathrm{rev}}+({r_{1}+r_{2}T})I/{A} \\
&+s\log\left[({t_{1}+t_{2}/T+t_{3}/T^{2}})I/{A} + 1\right], \label{eq:UI}
\end{split}
\end{equation}
\noindent
where $U^{\mathrm{stack}}$ and $U^{\mathrm{rev}}$ are stack voltage and the reversible voltage, respectively; $I$ and $T$ represent electrolytic current and stack temperature; $A$ is the electrode area; $r_{1}$, $r_{2}$, $s$, $t_{1}$, $t_{2}$ and $t_{3}$ are constant coefficients.

The DC load of the electrolytic stack $P^{\mathrm{stack}}$, as well as the operating ranges for $I$ and $T$, are as follows:
\begin{align}
P^\mathrm{stack}&=U^\mathrm{stack}I,  \label{eq:PI} \\
b^\mathrm{On}I^\mathrm{min}\leq I\leq b^\mathrm{On}&I^\mathrm{max},\ T^\mathrm{min}\leq T\leq T^\mathrm{max},  \label{eq:I-T}
\end{align}
\noindent
where $b^\mathrm{On}$ represents the production state; $I^\mathrm{max}$ and $I^\mathrm{min}$ are current limits; $T^\mathrm{max}$ and $T^\mathrm{min}$ are temperature limits.

According to Faraday's law, hydrogen production flow $Y^\mathrm{H_2}$ of a stack composed of $N^\mathrm{cell}$ cells satisfies (\ref{eq:h2}), and the Faraday efficiency $\eta^\mathrm{F}$ and overall energy conversion efficiency $\eta^\mathrm{ele}$ are determined in (\ref{eq:falad}) and (\ref{eq:eff}):
\begin{align}
  Y^\mathrm{H_2} &= \eta^\mathrm{F} {N^\mathrm{cell}I}/{(2F)}\times2\times3600/1000, \label{eq:h2} \\
  \eta^\mathrm{F} &= (I/A)^2/({f_1+\left(I/A\right)^2}) \times f_2, \label{eq:falad} \\
  \eta^\mathrm{ele} &= Y^\mathrm{H_2}\mathrm{L}^{\mathrm{hv}}/P^\mathrm{stack}\label{eq:eff},
 \end{align}
\noindent
where $F$ is the Faraday constant; $f_{1}$ and $f_{2}$ are constant coefficients; $\mathrm{L}^{\mathrm{hv}}$ is the lower heating value of hydrogen.

The load power $P^\mathrm{ele}$ of an ELZ includes not only $P^\mathrm{stack}$, but also auxiliary consumption $P^\mathrm{aux}$ for pumps, cooling, and the controller, calculated as follows:
\begin{align}
  P^\mathrm{ele}=P^\mathrm{stack}+&P^\mathrm{aux}, \label{eq:Pele} \\
  P^{\mathrm{aux}}=\left(1-b^{\mathrm{Idle}}\right)P^{\mathrm{cool}}/&\eta^{\mathrm{cool}}+b^{\mathrm{By}}P^{\mathrm{By}}, \label{eq:Paux}
\end{align}
\noindent
where $b^{\mathrm{Idle}}$ and $b^{\mathrm{By}}$ represent the idle and standby states; $P^{\mathrm{cool}}$ and $\eta^{\mathrm{cool}}$ are the cooling heat flow and its efficiency.

\subsection{The Active and Reactive Power Coupling of the Electrolyzers}
\label{sec:real-reactive}

The reactive power load of the ELZ is determined by $I$ and $T$ \cite{li2023multi}. In this section, using the commonly used 24-pulse TR as an example, we derive its expression.

First, we divide the reactive load $Q$ into the phase shift component $Q_{s}$ and distortion component $Q_{d}$ \cite{ruuskanen2020power}, calculated as (\ref{eq:Q})-(\ref{eq:Qs-Qd}). Then, we establish the relation between the AC voltage $U^{\mathrm{ac}}$ and $U^{\mathrm{stack}}$ as in (\ref{eq:Uac-dc}) \cite{li2023multi} and AC current $I^{\mathrm{ac}}$ as in (\ref{eq:Pac-dc}). The loss of the rectifier $P^{\mathrm{loss}}$ are approximated as (\ref{eq:Ploss}) \cite{zheng2020discrete}. Subsequently, the expression for $Q(U^{\mathrm{ac}},I,T)$ can be obtained, and the detailed derivation is provided in the Appendix.
\begin{align}
  Q &=\sqrt{Q_{s}^{2}+Q_{d}^{2}}, \label{eq:Q} \\
  Q_{s}=\sqrt{3}U^{\mathrm{ac}}I^{\mathrm{ac}}&\sin\varphi,\ Q_{d}=\sqrt{3}U^{\mathrm{ac}}I^{\mathrm{ac}}{\sqrt{1-\nu^{2}}}/{\nu}, \label{eq:Qs-Qd} \\
  U^{\mathrm{stack}} & =U^{\mathrm{ac}}{2.44}/{K} \times \cos\varphi, \label{eq:Uac-dc} \\
  I^{\mathrm{ac}} & =\dfrac{U^\mathrm{stack}I+P^\mathrm{loss}}{\sqrt{3}U^\mathrm{ac}\cos\varphi}, \label{eq:Pac-dc} \\
  P^{\mathrm{loss}} & =a_2^{\mathrm{loss}}I^2+a_1^{\mathrm{loss}}I+a_0^{\mathrm{loss}}, \label{eq:Ploss}
\centering
\end{align}
\noindent
where $\varphi$ is the firing angle; $\nu$ is the harmonic factor corresponding to 24-TR; $K$ is the turns ratio of the transformer; $a_0^{\mathrm{loss}}$, $a_1^{\mathrm{loss}}$, and $a_2^{\mathrm{loss}}$ are the fitting coefficients for $P^{\mathrm{loss}}$.

\subsection{Example of Trade-offs Between Energy Conversion Efficiency and Reactive Power in Scheduling Multiple Electrolyzers}
\label{sec:conflict}

\begin{figure}[tb]
  \centering
  \includegraphics[width=3.3in]{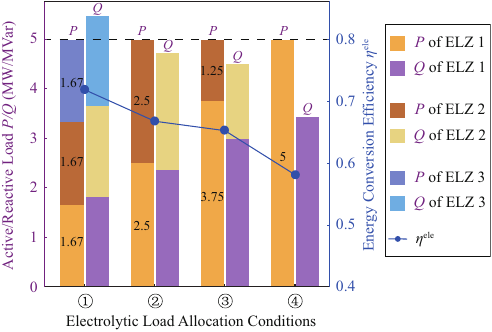} \vspace{-2.5pt}
  \caption{Reactive power $Q$ and power-to-hydrogen energy conversion efficiency $\eta^{\mathrm{ele}}$ under 4 conditions of electrolytic load allocation among 3 electrolyzers.}
  \label{fig:conflict}\vspace{-10pt}
\end{figure}

To provide a clear illustration of how load allocation among multiple ELZs affects energy conversion efficiency and reactive power, we present an example. Given a power supply of 5 MW, we consider 4 operating conditions ($T=80$ \textcelsius): \textcircled{1} 3 ELZs operating at 1.67 MW each; \textcircled{2} 2 ELZs at 2.5 MW each; \textcircled{3} 1 ELZ at 3.75 MW and 1 ELZ at 1.25 MW; \textcircled{4} 1 ELZ at 5 MW. The efficiency $\eta^\mathrm{ele}$ and reactive load $Q$ are shown in Fig. \ref{fig:conflict}.

Clearly, $\eta^\mathrm{ele}$ and $Q$ exhibit significant variations as load allocation among the ELZs varies, even when the total active load remains constant. When the power is evenly allocated, $\eta^\mathrm{ele}$ is higher, but $Q$ also increases. In contrast, concentrating all the power on a single ELZ results in the lowest values for both $\eta^\mathrm{ele}$ and $Q$. As a result, when only aims at maximizing $\eta^\mathrm{ele}$ as \cite{varela2021modeling,QIU2023EXTEND}, excessive $Q$ may lead to voltage drops and increased network losses. In contrast, only optimizing the reactive power flow may lead to low energy conversion efficiency.

\subsection{Dynamic Thermal Model}

The temperature dynamics of the ELZ can be described by a first-order model \cite{ulleberg2003modeling}, as:
\begin{align}
  C_m^{\mathrm{ele}}({T_{m,t+1}-T_{m,t}})/{\Delta t} &=P_{m,t}^{\mathrm{gen}}-P_{m,t}^{\mathrm{diss}}-P_{m,t}^{\mathrm{cool}}, \label{eq:T-P}
\end{align}
\noindent
where the electrolytic heat $P_{m,t}^{\mathrm{gen}}$, heat dissipation $P_{m,t}^{\mathrm{diss}}$ and cooling heat $P_{m,t}^{cool}$ are constrained by
\begin{align}
  \begin{cases}
    P_{m,t}^{\mathrm{gen}} =I_{m,t}\left(U_{m,t}^{\mathrm{stack}}-U^{\mathrm{th}}\right), \\
    P_{m,t}^{\mathrm{diss}}={(T_{m,t}-T^{\mathrm{am}})}/{R_{m}^{\mathrm{diss}}},  \\
    0\leq P_{m,t}^{cool} \leq C^{cool}\left(T_{m,t}-T^{cool}\right), \label{eq:Pcool}
  \end{cases}
\end{align}
\noindent
where $C_m^{\mathrm{ele}}$ and $R_{m}^{\mathrm{diss}}$ are heat capacity and thermal resistance of the \textit{m}th ELZ; $U^{\mathrm{th}}$ is the thermal neutral voltage; $T^{\mathrm{am}}$ is the ambient temperature; $C^{cool}$ and $T^{cool}$ are heat capacity and temperature of cooling water.

\subsection{Operational State Transitions}

The operating states of the ELZ include Production, Standby, and Idle, denoted as $b_{m,t}^{\mathrm{On}}$, $b_{m,t}^{\mathrm{By}}$, and  $b_{m,t}^{\mathrm{Idle}}$, respectively, which satisfy the following logic expressions, as
\begin{align}
  \begin{cases}
    b_{m,t}^{\mathrm{On}}+b_{m,t}^{\mathrm{By}}+b_{m,t}^{\mathrm{Idle}}=1,  \\
    b_{m,t}^{\mathrm{On}}+b_{m,t}^{\mathrm{By}}+b_{m,t-1}^{\mathrm{Idle}}-1\leq b_{m,t}^{\mathrm{SU}},  \\
    b_{m,t-1}^{\mathrm{On}}+b_{m,t-1}^{\mathrm{By}}+b_{m,t}^{\mathrm{Idle}}-1\leq b_{m,t}^{\mathrm{SD}},  \\
    -b_{m,t-2}^{\text{Idle}}-b_{m,t}^{\text{Idle}}+b_{m,t-1}^{\text{Idle}} \leq 0, \label{eq:bbb}
  \end{cases}
\end{align}
where $b_{m,t}^{\mathrm{SU}}$ and $b_{m,t}^{\mathrm{SD}}$ represent startup and shutdown actions, and more detailed description can be seen in \cite{varela2021modeling} or \cite{QIU2023EXTEND}.

\section{Coordinated Active-Reactive Power Management for the ReP2H System}
\label{sec:management}

Aiming at maximizing the profit of hydrogen production via coordinated active and reactive power management of the ELZs and other resources, the detailed model is elaborated below.

\subsection{Network Power Flow of the ReP2H System }
\label{sec:network}

The network power flow constraints of the ReP2H system are described by the classical DistFlow model (see \cite{chai2019hierarchical}), where power injections at each bus are determined by:
\begin{align}
& \begin{cases}
\text{DistFlow Equations \cite{chai2019hierarchical}}, \\
p_{j,t}=P_{j,t}^\mathrm{WT}+P_{j,t}^\mathrm{PV}+P_{j,t}^\mathrm{ES,out}-P_{j,t}^\mathrm{ES,in}-p_{j,t}^\mathrm{L},\\
q_{j,t}=Q_{j,t}^\mathrm{WT}+Q_{j,t}^\mathrm{PV}+Q_{j,t}^\mathrm{ES}+Q_{j,t}^\mathrm{CB}+Q_{j,t}^\mathrm{C}-q_{j,t}^\mathrm{L},\\
p_{j,t}^\mathrm{L}=\sum_{ m\in\mathbb{M}}\left(P_{m,t}^\mathrm{ele}+P_{m,t}^\mathrm{loss}\right),\\
q_{j,t}^\mathrm{L}=\sum_{ m\in\mathbb{M}}Q, \label{eq:node pq}
\end{cases}
\end{align}
where $P_{j,t}^\mathrm{WT}$, $Q_{j,t}^\mathrm{WT}$, $P_{j,t}^\mathrm{PV}$, and $Q_{j,t}^\mathrm{PV}$ represent the active and reactive power of WT and PV at bus j, respectively; $P_{j,t}^\mathrm{ES,in}$, $P_{j,t}^\mathrm{ES,out}$, and $Q_{j,t}^\mathrm{ES}$ are the charging, discharging, and reactive power of ES; $Q_{j,t}^\mathrm{CB}$ and $Q_{j,t}^\mathrm{C}$ represent the reactive power of the CBs and SVC; $p_{j,t}^\mathrm{L}$ and $q_{j,t}^\mathrm{L}$ are the active and reactive loads of the hydrogen plant; $\mathbb{M}$ is the set of ELZs.

\subsection{Reactive Power Resources}
\label{sec:resources}

The controllable reactive power sources, including WTs, PV plants, ESs, capacitor banks (CBs), and the static var compensator (SVC), are modeled as follows.

\subsubsection{Reactive Power Capacity of WTs \cite{zhang2016reactive}}

\begin{align}
  \begin{cases}
    Q_{j,t}^\mathrm{WT} \le -0.58P_{j,t}^\mathrm{WT} + 0.91S_{j}^\mathrm{WT}, \\
    Q_{j,t}^\mathrm{WT} \ge 1.24P_{j,t}^\mathrm{WT} - 0.91S_{j}^\mathrm{WT}, \\
    0 \le P_{j,t}^\mathrm{WT}\le S_{j}^\mathrm{WT},
  \end{cases}
\end{align}
\noindent
where $S_{j}^\mathrm{WT}$ is the installation capacity of the $j$th \upshape{WT}.

\subsubsection{Reactive Power Capacity of PV plants \cite{chai2019hierarchical}}

\begin{align}\begin{cases}
  \big\Vert \big(P_{j,t}^{\mathrm{PV}} ,  Q_{j,t}^{\mathrm{PV}} \big) \big\Vert_2 \le S_{j}^{\mathrm{PV}},\\
  -P_{j,t}^{\mathrm{PV}}\tan\theta\le Q_{j,t}^{\mathrm{PV}}\le P_{j,t}^{\mathrm{PV}}\tan\theta,
\end{cases}\end{align}
\noindent
where $S_{j}^\mathrm{PV}$ is the installation capacity of the $j$th PV plant; $\theta$ is the maximum power factor angle.

\subsubsection{ES Operational Constraints}

\begin{align}\hspace{-6pt}
  \begin{cases}
    \big\Vert \big( P_{j,t}^{\mathrm{ES,in}} , Q_{j,t}^{\mathrm{ES}} \big) \big\Vert_2 \le S_{j}^{\mathrm{ES}}, \\
    \Vert \big(P_{j,t}^{\mathrm{ES,out}} , Q_{j,t}^{\mathrm{ES}}\big) \big\Vert_2 \le S_{j}^{\mathrm{ES}}, \\
    0\leq P_{j,t}^{\mathrm{ES,in}}\leq b_{j,t}^{\mathrm{ES,in}}P_{\mathrm{max}}^{\mathrm{ES,in}}, \\
    0\leq P_{j,t}^{\mathrm{ES,out}}\leq b_{j,t}^{\mathrm{ES,out}}P_{\mathrm{max}}^{\mathrm{ES,out}}, \\
    b_{j,t}^{\mathrm{ES,in}}+b_{j,t}^{\mathrm{ES,out}} \leq 1 ,\\
    SOC_{j,t}=SOC_{j,t-1}+\big(\eta^{\mathrm{ES,in}}P_{j,t}^{\mathrm{ES,in}}-{P_{j,t}^{\mathrm{ES,out}}}/{\eta^{\mathrm{ES,out}}}\big)\Delta t, \\
    \underline{SOC} \le SOC_{j,t} \le \overline{SOC},
  \end{cases}\hspace{-36pt}
\end{align}
\noindent
where $S_{j}^{\mathrm{ES}}$ is the convertor capacity; $P_{\mathrm{max}}^{\mathrm{ES,in}}$ and $P_{\mathrm{max}}^{\mathrm{ES,out}}$ are the charging and discharging power limits; $b_{j,t}^{\mathrm{ES,in}}$ and $b_{j,t}^{\mathrm{ES,out}}$ depict the operation state; $\eta^{\mathrm{ES,in}}$ and $\eta^{\mathrm{ES,out}}$ are efficiencies; $SOC_{t}$, $\underline{SOC}$, and $\overline{SOC}$ are the state of charge and its limits.

\subsubsection{CB and SVC Operation Constraints}
\begin{align}
  &\begin{cases}Q_{j,t}^{\mathrm{CB}}=n_{j,t}^{\mathrm{CB}}\upsilon_{j,t}\Delta Q_{j}^{\mathrm{CB}}, \\
    0\le n_{j,t}^{\mathrm{CB}}\le\overline{n}_{j}^{\mathrm{CB}}, \ \left|n_{j,t+1}^{\mathrm{CB}}-n_{j,t}^{\mathrm{CB}}\right|\le n_{j}^{\mathrm{CB,max}}, \label{eq:CB}
  \end{cases} \\
  &-Q_{j,\mathrm{max}}^{\mathrm{C}}\leq Q_{j,t}^{\mathrm{C}}\leq Q_{j,\mathrm{max}}^{\mathrm{C}}, \label{eq:svc}
\end{align}
\noindent
where $\upsilon_{j,t}$ is the squared voltage magnitude; $\Delta Q_{j}^{\mathrm{CB}}$ is the capacity of a single CB at rated voltage; $n_{j,t}^{\mathrm{CB}}$ and $\overline{n}_{j}^{\mathrm{CB}}$ represent the number of CBs in operation and its limit; $n_{j}^{\mathrm{CB,max}}$ is the maximal number of CBs that can be switched in a single operation; $Q_{j,t}^{\mathrm{C}}$ and $Q_{j,\mathrm{max}}^{\mathrm{C}}$ are the reactive power and limits of the SVC.

\subsection{The Overall Scheduling Problem }
\label{sec:problem}

The objective is maximizing the total profit $R$, which includes hydrogen production revenue, startup and shutdown costs, and CB action cost, formulated by
\begin{align}
R=\sum_{m\in{\mathbb{M}}}\sum_{t\in\mathbb{T}}\left\{\begin{matrix}
c^{\mathrm{H}_2}Y_{m,t}^{\mathrm{H}_2}-c^{\mathrm{SU}}b_{m,t}^{\mathrm{SU}}-c^{\mathrm{SD}}b_{m,t}^{\mathrm{SD}}\\
-c^{\mathrm{CB}}{\sum}_j\Big|n_{j,t+1}^{\mathrm{CB}}-n_{j,t}^{\mathrm{CB}}\Big|
\end{matrix}\right\}, \label{eq:R}
\end{align}
\noindent
where $\mathbb{T}$ is the scheduling horizon; $c^{\mathrm{H}_2}$, $c^{\mathrm{SU}}$, $c^{\mathrm{SD}}$ and $c^{\mathrm{CB}}$ are the hydrogen price, costs for startup and shutdown of ELZs and switching of CBs, respectively.

Overall, the decision variables include $\bm{b}$, $\bm{I}$, $\bm{T}$, $\bm{Q}^\mathrm{WT}$, $\bm{Q}^\mathrm{PV}$, $\bm{Q}^\mathrm{ES}$, $\bm{P}^\mathrm{ES,in}$, $\boldsymbol{P}^\mathrm{ES,out}$, $\bm{n}^\mathrm{CB}$ and $\bm{Q}^\mathrm{C}$, summarized as $\boldsymbol{x}$, and the proposed scheduling problem is formulated as follows:
\begin{align}
  \max\limits_{\bm{x}}\;(\ref{eq:R}),\ \   \text{s.t.} \;\;(\ref{eq:UI})-(\ref{eq:svc}). \label{eq:problem}
\end{align}

Due to the nonlinearity of $P^\mathrm{stack}$, and $Q$, etc. with respect to $I$, (\ref{eq:problem}) is a mixed-integer nonlinear programming. To solve it, we transform it into a mixed-integer second-order cone programming (MISOCP) via polynomial approximation, piecewise linearization, and the big-M method.

\section{Case Studies}
\label{sec:case}

\subsection{Case Settings}
\label{sec:setting}

\begin{figure}[tb]
  \centering
  \includegraphics[width=3.4 in]{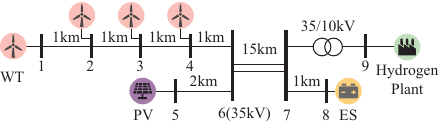}\vspace{-4.5pt}
  \caption{The topology of the off-grid ReP2H system used for case study.}
  \label{fig:case}\vspace{-2pt}
\end{figure}

An off-grid ReP2H system shown in Fig. \ref{fig:case} is used for the case studies. It comprises 4 $\times$ 6.25 MW WTs connected to buses 1-4, 1 $\times$ 5 MW PV plant at bus 5, 1 $\times$ 2.5 MW/5 MWh ES at bus 8, 6 $\times$ 0.5 MVar CBs and 1 $\times$ 1 MVar SVC at bus 7, and 4 $\times$ 5 MW alkaline ELZs at bus 9. The case parameters are shown in Table \ref{tab:para}. Simulations are performed on \emph{Wolfram Mathematica 12.3}, and the MISOCP is solved by \emph{Gurobi 9.5.2}.

\begin{table}[tb]\scriptsize
  \renewcommand{\arraystretch}{1.4}
  \caption{Parameter Setting of the Case Studies}\vspace{-6pt}
  \label{tab:para}
  \centering
  \begin{tabular}{cccc}
  \hline\hline
  Parameter                               & Symbol                & Value        \\
  \hline
  Scheduling horizon and step length      & $\mathbb{T}$ / $\Delta t$               & 24 / 1 h                    \\
  Current limit of the electrolysis stack & $I^\mathrm{max}$ / $I^\mathrm{min}$     & 12 / 3 kA                 \\
  ELZ temperature limits                  & $T^\mathrm{max}$ / $T^\mathrm{min}$     & 80 / 25 \textcelsius    \\
  Hydrogen price                          & $c^{\mathrm{H}_2}$                      & 29 CNY/kg  \cite{varela2021modeling}              \\
  ELZ startup and shutdown costs          & $c^{\mathrm{SU}}$ / $c^{\mathrm{SD}}$   & 1000 / 0 CNY \cite{varela2021modeling}           \\
  CB action cost                          & $c^{\mathrm{CB}}$                       & 2 CNY  \cite{chai2019hierarchical}                 \\
  Bus voltage limits                      & $\overline{U}_i$ /$\underline{U}_i$     & 1.05 / 0.95  p.u.               \\
  \hline\hline
  \end{tabular}\vspace{-14pt}
\end{table}

\begin{figure}[tb]
  \centering
  \includegraphics[width=3.47in]{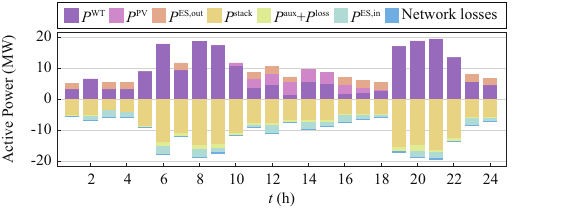}\vspace{-8pt}
  \caption{Active power supply and demand in the ReP2H system.}\vspace{-3pt}
  \label{fig:balance}
  \vspace{10pt}
  \includegraphics[width=3.47in]{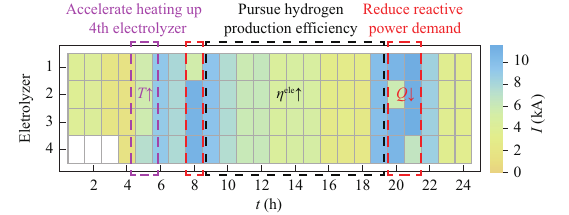}\vspace{-8pt}
  \caption{Operational states and active load allocation of the electrolyzers.}
  \label{fig:distribution}
  \vspace{10pt}
  \centering
  \includegraphics[width=3.47in]{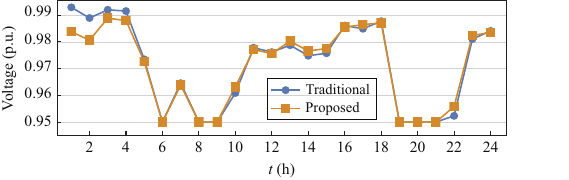}\vspace{-8pt}
  \caption{Voltage magnitude at the load bus.}
  \label{fig:voltage}
  \vspace{10pt}
  \includegraphics[width=3.47in]{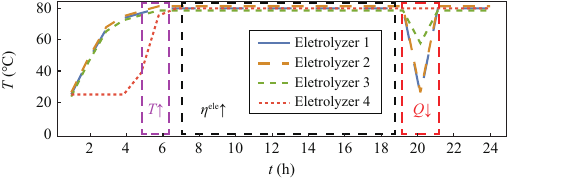}\vspace{-8pt}
  \caption{Stack temperatures of the electrolyzers.}
  \label{fig:temperature}
  \vspace{10pt}
  \centering
  \includegraphics[width=3.47in]{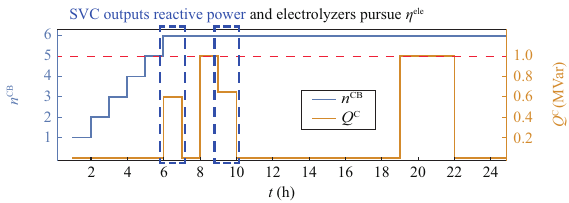}\vspace{-8pt}
  \caption{States of the reactive power compensators.}
  \label{fig:compensation}\vspace{-12pt}
\end{figure}

\subsection{Simulation Result of the Proposed Coordinated Power Energy Management Method}
\label{sec:analysis}

With the PV and wind power generation given in Fig. \ref{fig:balance}, the optimized on-off switching sequence and load allocation of the 4 ELZs by the proposed method are presented in Fig. \ref{fig:distribution}. The active power balance of the system is shown in Fig. \ref{fig:balance}, with voltage profile at the load bus shown in Fig. \ref{fig:voltage}. Moreover, Fig. \ref{fig:temperature} displays the temperatures of the ELZs.

As shown in Fig. \ref{fig:distribution}, most of the time, the active load is evenly allocated to the ELZs to optimize energy conversion efficiency $\eta^\mathrm{ele}$. However, at $t=\{8,20,21\}$, the load is concentrated in 2 or 3 ELZs to reduce total reactive demand $Q$ and avoid voltage violation. The detailed explanation can be seen in Section \ref{sec:conflict}. Due to the high startup cost of the ELZs, it is generally not economically feasible to reduce $Q$ by shutting down ELZs. From Fig. \ref{fig:balance} and Fig. \ref{fig:voltage}, we can see that if the hydrogen plant focuses on improving $\eta^\mathrm{ele}$ at these moments, it may lead to an excessive $Q$ and increased network losses.

Moreover, at $t=5$, we can observe the temperature of the 4th ELZ is lower than the others due to its late startup. Therefore, the proposed scheduling approach commanded a full load to speed up heating up to increase $\eta^\mathrm{ele}$. Also, the electrolyzers reduces temperature at $t=\{20,21\}$ to reduce reactive power demand, which accords with the negative correlation between $T$ and $Q(U^{\mathrm{ac}},I,T)$ derived in Section \ref{sec:real-reactive}.

Fig. \ref{fig:compensation} shows the outputs of the var compensators. Comparing $t=8$ and $t=9$, it is observed that when the reactive power is sufficient, the hydrogen plant will not decrease $Q$ at $t = 9$, in contrast to the insufficient condition at $t = 8$. In other words, we know that the benefits of reducing network losses cannot offset the losses caused by reducing $\eta^\mathrm{ele}$ here.

\subsection{Comparison Between the Proposed Method and the Traditional Scheduling Method}
\label{sec:comparison}

\begin{figure}[tb]
  \centering
  \includegraphics[width=3.47in]{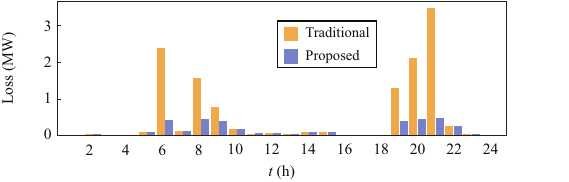}\vspace{-8pt}
  \caption{Comparison of network losses between the proposed and the traditional scheduling methods.}
  \label{fig:loss}\vspace{-12pt}
\end{figure}

Table \ref{tab:comparision} compares the proposed method and the traditional hydrogen plant scheduling method \cite{QIU2023EXTEND}. The network losses is specifically shown in Fig. \ref{fig:loss}. As can be seen, when the nodal voltage is at a critical value (see Fig. \ref{fig:voltage}), the proposed method exhibits a significant reduction in the network losses. This validates that the proposed method better adapts to various operating conditions and improves the profitability.

\begin{table}[tb]\scriptsize
  \renewcommand{\arraystretch}{1.4}
  \caption{Comparison on Hydrogen Production, Network Losses, and Overall Profit Between the Proposed and the Traditional Scheduling Methods Under the Base Case}\vspace{-6pt}
  \label{tab:comparision}
  \centering
  \begin{tabular}{cccc}
  \hline\hline
  Method                          & Hydrogen Production       &Network Loss Ratio  & Profit    \\
  \hline
  Traditional              &  3974.61 kg              & 4.32\%      & 111264 CNY   \\
  Proposed               & 4153.54 kg               & 1.31\%        & 115818 CNY  \\ \hline
  Comparison                        & +4.5\%                 & -3.01\%   & +4.09\%     \\
  \hline\hline
  \end{tabular}\vspace{-3pt}
\end{table}

\subsection{Comparison Under Various Renewable Power Scenarios}
\label{sec:scenarios}

\begin{figure}[tb]
  \centering\vspace{-8pt}
  \includegraphics[width=3.43in]{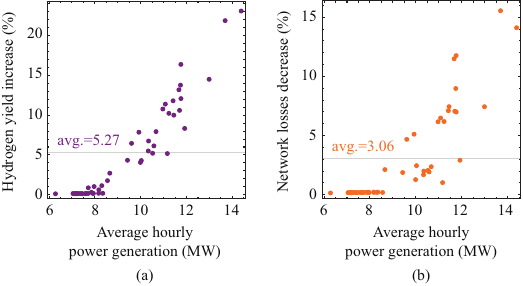}\vspace{-5pt}
  \caption{Improvement of the proposed coordinated active-reactive power management method under various scenarios of renewable power in terms of (a) hydrogen production (b) reduction of network losses.}
  \label{fig:50}
\end{figure}

Finally, the proposed method is validated under 50 scenarios of wind and solar output. The data come from Inner Mongolia, China.
As shown in Fig. \ref{fig:50}(a) and \ref{fig:50}(b), hydrogen production is on average increased by 5.27\%, and the network losses is reduced by 3.06\% compared to the traditional method. Particulary, the improvement is more significant during periods with high wind and solar generation, when the system is more stressed in both active and reactive power flow.

\section{Conclusions}
\label{sec:conclusion}

This study reveals a trade-off between P2H energy conversion efficiency and network losses in ReP2H systems caused by the reactive power demands of electrolyzers. To address this, a coordinated active-reactive power management method is proposed. Case studies show that the proposed method significantly improves hydrogen yield and overall profitability compared to traditional scheduling methods only optimizing active power.

Future study may focus on dealing with uncertainty from the renewable power sources and scalability for applications in large-scale systems composed of more electrolyzers.

\section{Acknowledgement}

Financial support from National Key R\&D Program of China (2021YFB4000503), National Natural Science Foundation
of China (52377116 and 52207116), and China Postdoctoral Science Foundation (2022M711758) is gratefully acknowledged.

\appendix
\label{sec:appendix}
\setcounter{equation}{0}
\renewcommand{\theequation}{\ref{sec:appendix}\arabic{equation}}

Combining (\ref{eq:Q}) and (\ref{eq:Qs-Qd}), we have
\begin{align}
  Q= \sqrt{3} U^\text{ac}I^\text{ac}\sqrt{\sin^2\varphi + ({1-\nu^2})/{\nu^2}}. \label{eq:A1}
\end{align}

Then, substituting (\ref{eq:Pac-dc}) into (\ref{eq:A1}) yields
\begin{align}
  Q = ({U^\text{stack} I + P^{\text{loss}}} )/{\cos\varphi} \times \sqrt{\sin^2\varphi+ ({1-\nu^2})/{\nu^2}}. \label{eq:A2}
\end{align}

Also, using (\ref{eq:Uac-dc}) we have
\begin{equation}
 \cos\varphi=\frac{K}{2.44}\frac{U^{\mathrm{stack}}}{U^{\mathrm{ac}}}\Leftrightarrow\varphi=\arccos\left(\frac{K}{2.44}\frac{U^{\mathrm{stack}}}{U^{\mathrm{ac}}}\right). \label{eq:A3}
\end{equation}

By substituting equations (\ref{eq:Ploss}) and (\ref{eq:A3}) into (\ref{eq:A2}), we obtain the expression of $Q$ in terms of $U^{\mathrm{stack}}$, $I$, and $U^{\mathrm{ac}}$, as
\begin{align}
  Q = & {2.44\Big(U^\text{stack}I+a_2^{1\text{oss}}I^2+a_1^{1\text{oss}}I+a_0^{1\text{oss}}\Big)U^\text{ac}}/{( K U^\text{stack})} \nonumber \\
  \times & {\sqrt{   \sin^ 2 \big[ \text{arccos} \big({KU^{\text{stack}}}/({2.44U^{\text{ac}}}) \big) \big]+ ({1-\nu^2})/{\nu^2}}}. \label{eq:A4}
\end{align}

Substituting (\ref{eq:UI}) into (\ref{eq:A4}), we obtain $Q(U^{\mathrm{ac}},I,T)$.

\bibliographystyle{IEEEtran}
\bibliography{IEEEabrv,ReP2H-ORPF-GM}

\end{document}